\documentclass[10pt]{amsart}
\usepackage{amsmath,amssymb,amscd}
\usepackage[all]{xy}

\newtheorem{Df}{Definition}[section]
\newtheorem{Te}[Df]{Theorem}
\newtheorem{Po}[Df]{Proposition}
\newtheorem{Cr}[Df]{Corollary}
\newtheorem{Lm}[Df]{Lemma}
\newtheorem{Ca}[Df]{Claim}
\newtheorem{Cn}[Df]{Conjecture}

\newcommand{\Bdf}{\begin{Df}}
\newcommand{\Edf}{\end{Df}}
\newcommand{\Bte}{\begin{Te}}
\newcommand{\Ete}{\end{Te}}
\newcommand{\Bpo}{\begin{Po}}
\newcommand{\Epo}{\end{Po}}
\newcommand{\Bcr}{\begin{Cr}}
\newcommand{\Ecr}{\end{Cr}}
\newcommand{\Blm}{\begin{Lm}}
\newcommand{\Elm}{\end{Lm}}
\newcommand{\Bca}{\begin{Ca}}
\newcommand{\Eca}{\end{Ca}}
\newcommand{\Bcn}{\begin{Cn}}
\newcommand{\Ecn}{\end{Cn}}
\newcommand{\Bdm}{{\it Proof.}\ }
\newcommand{\Rm}{{\it Remark \arabic{section}.\arabic{Df} \ }}

\begin{document}
\title{Gerasimov's theorem and $N$-Koszul algebras}

\author{Roland Berger}
\date{}

\maketitle
\thispagestyle{empty}
\begin{abstract}
This paper is devoted to graded algebras $A$ having a single homogeneous relation. 
We give a criterion for $A$ to be $N$-Koszul where $N$ is the degree of the relation. 
This criterion uses a theorem due to Gerasimov. As a consequence of the criterion, 
some new examples of $N$-Koszul algebras are presented. We give an alternative proof of Gerasimov's theorem for $N=2$, 
which is related to Dubois-Violette's theorem concerning a matrix description of the Koszul and AS-Gorenstein algebras of global dimension 2. We determine which of the PBW deformations of a symplectic form are Calabi-Yau.
\end{abstract}

\section{Introduction}
The first elementary connection between Koszul duality and combinatorics is probably the computation 
of the Hilbert series of a polynomial algebra as the inverse of the Hilbert series of the corresponding 
Grassmann algebra. Actually this numerical Hilbert series combinatorics holds for any Koszul algebra, 
including the $N$-generalisation of Koszul algebras (see e.g.~\cite{dvp:plac, kriegk:cras}). 
Recently, Phung Ho Hai and Martin Lorenz showed that the MacMahon Master Theorem (MMT) is a consequence 
of another type of Hilbert series combinatorics~\cite{phhml:2MMT}, which is based on the Grothendieck ring of 
comodules over Manin's bialgebra 
of any 2-Koszul algebra (specialised to a polynomial algebra if we want to recover the original MMT). 
This new type of Hilbert series combinatorics can 
be called comodule Hilbert series combinatorics. Comodule Hilbert series combinatorics reduces to numerical Hilbert series combinatorics by taking the dimension of comodules (comodules are assumed to be finite dimensional as vector spaces). More recently, Phung Ho Hai, Benoit Kriegk and Martin Lorenz showed that comodule 
Hilbert series combinatorics holds 
for any $N$-Koszul algebra (or superalgebra as well)~\cite{phhml:NMMT}, including certain explicit combinatorial applications 
(see also~\cite{ep:algMMT}). 

Another interesting combinatorics linked to $N$-Koszul algebras can be traced back to Backelin's PhD thesis~\cite{bafr:kalg} 
and can be called 
combinatorics of $V$'s and $R$'s. Here $V$ and $R$ denote respectively the space of generators and 
the space of relations of the algebra $A$, 
where $A$ is assumed to be $N$-homogeneous (i.e., $R$ is a subspace of $V^{\otimes N}$). Then $A$ is Koszul if and only if certain 
(infinitely many) triples of vector spaces are distributive and, when $N>2$, an extra condition holds 
(see~\cite{bafr:kalg} if $N=2$, and~\cite{rb:nonquad} for any $N\geq 2$). 
Each vector space in such a triple is either a sum or an intersection of tensor products of $V$'s and $R$'s with the same 
number of factors, $R$ occurring only once in each tensor product. Formulas giving the sums or the intersections are very specific. 
Recall that a triple $(E,F,G)$ of subspaces of a vector 
space is said to be \emph{distributive} if the relation $E\cap (F+G)=(E\cap F)+(E\cap G)$ holds. 

Gerasimov's theorem is a result in the combinatorics of $V$'s and $R$'s. The algebra $A$ is still $N$-homogeneous, 
but now $R$ is one-dimensional. 
Then Gerasimov's theorem asserts that $A$ is \emph{distributive}, meaning that 
\emph{any} triple of sums or intersections of tensor products of $V$'s and $R$'s (with the same 
number of factors in the tensor products, $R$ occurring only once) is distributive~\cite{gera:distrib} (in the present text, 
all the generators have degree 1, but let us notice that 
Gerasimov's theorem holds for generators with positive degrees which may be unequal and distinct from 1).
Using the previous characterisation of Koszulity, we get the following criterion, presented in Section 2. 
\Bte \label{crit}
Let $V$ be a vector space over a field $k$ and let $N$ be an integer $\geq 2$. Let $R$ be a one-dimensional subspace of
$V^{\otimes N}$. Then the graded algebra $A$ defined by the space $V$ of generators (all assumed to be of degree 1) and the space
$R$ of relations is $N$-Koszul if and only if we have
\begin{equation}
(R\otimes V^{\otimes m}) \cap (V^{\otimes m}\otimes R) \subseteq V^{\otimes (m-1)}\otimes R \otimes
V, \ \mathrm{for} \ m=2, \ldots , N-1.
\end{equation} 
In particular, $A$ is Koszul when $N=2$.
\Ete

When $N=2$, this result is the noncommutative version of the following well-known result:
every commutative graded algebra whose generators are of degree 1, subject to a single relation which is homogeneous quadratic, 
is Koszul~\cite{bafr:kalg}.

In Section 3, we will illustrate Theorem \ref{crit} by the following.
\Bcr \label{anti}
We keep the notations and assumptions of \emph{Theorem \ref{crit}}. Assume moreover that $k$ has characteristic zero, $V$ is 
finite-dimensional and $R$ is generated by an antisymmetric tensor. Then $A$ is $N$-Koszul and the global dimension of $A$ is 2.
\Ecr

When the antisymmetric tensor is an antisymmetriser, another proof of the Koszulity in this corollary is presented in Section 3. It is based on a basic result 
(interesting for its own sake) concerning  
free products of distributive and $N$-Koszul algebras. However, although this proof does not use Gerasimov's theorem, it 
uses the fact that the antisymmetriser algebras (also called $N$-symmetric algebras~\cite{phhml:NMMT}) are $N$-Koszul and distributive, and this fact is not 
obvious either~\cite{rb:nonquad}. 

A criterion analogous to that of Theorem \ref{crit} was already known for monomial algebras~\cite{rb:nonquad}. The reason is the same 
as in the one-dimensional situation, namely any monomial algebra is distributive (but here the proof of the distributivity 
is immediate, unlike the proof of Gerasimov's theorem). So only the extra condition is required 
in both situations in order to have an $N$-Koszul algebra. Combining both situations, we will obtain in Section 4 the following 
corollary.
\Bcr
We keep the notations and assumptions of \emph{Theorem \ref{crit}}, with $V$ finite-dimensional. Let $(x_1, \ldots, x_n)$ be a basis of $V$, $n\geq 1$. 
Assume that $R$ is generated by a monomial $f=x_{i_1} \ldots x_{i_N}$, $N\geq 2$. Then $A$ is $N$-Koszul if and only if there exists no $m$ 
in $\{2,\ldots , N-1\}$ such that 
$$f=(x_{i_1} \ldots x_{i_m})^q\, x_{i_1} \ldots x_{i_r},$$
where $N=mq+r$ is the division of $N$ by $m$ with remainder $r$, and where $i_1,\ldots ,i_m$ are not all equal.
\Ecr

Gerasimov's theorem is considered an important contribution to combinatorial ring theory. For example, as shown by 
Backelin~\cite{bac:ratio}, the rationality of the Hilbert series of algebras with a single relation 
is a striking consequence of Gerasimov's theorem. 
The proof given by Gerasimov in~\cite{gera:distrib} of his theorem is rather long (36 pages) and is situated 
at the frontiers of algebra, logic and combinatorics, using an \emph{ad hoc} formal calculus with rules of inference and involving 
provability and decidability. Even if the interactions 
with logic seem exciting, it would be satisfactory to have an algebraic proof of Gerasimov's theorem, whose statement is purely 
algebraic. In Section 5, we will provide such an algebraic proof of Gerasimov's theorem in the quadratic case, i.e., for $N=2$.
The point is that distributivity for any 2-homogeneous graded algebra is equivalent to Koszulity~\cite{bafr:kalg}, so that we can use homological algebra when $N=2$, and this gives a very different proof. 

More precisely, it suffices to prove that a graded algebra $A$ with
generators $x_1, \ldots, x_n$ subject to the single relation $f=0$, where 
$$f=\sum_{1\leq i,j \leq n} f_{ij}\, x_i x_j$$
is Koszul. It turns out that Michel Dubois-Violette has investigated such graded algebras $A$ using homological algebra, 
showing that $A$ is Koszul and AS-Gorenstein of global dimension 2 if the matrix $M=(f_{ij})$ of the 
coefficients of $f$ is invertible~\cite{mdv:multi}. In Section 5, we make precise the proof of Koszulity 
and extend it to some noninvertible matrices $M$. Next we use Bergman's confluence for the remaining $M$'s. 
Global dimension, Hilbert series and Gelfand-Kirillov dimension are computed systematically, 
and we recover the fact that $A$ is AS-Gorenstein if and only if $M$ is invertible.

Finally, Section 6 is devoted to Poincar\'e-Birkhoff-Witt (PBW) deformations of various $N$-Koszul graded algebras appearing in the previous sections. We shall determine in Theorem \ref{cypbw} which of the PBW deformations of a symplectic form are Calabi-Yau, the symplectic form being viewed as a homogeneous quadratic relation. This gives examples of graded Calabi-Yau algebras with some PBW deformations that are not Calabi-Yau. In general the Calabi-Yau property is preserved under PBW deformations when only constants are added to relations (Theorem \ref{cyconst}).
\\ \\
\textbf{Acknowledgements.} I am very grateful to Dmitriy Rumynin for pointing out Gerasimov's paper. I also thank Jacques Alev for 
suggesting to use free products of algebras for Corollary \ref{anti}, and Julien Bichon for pointing out reference~\cite{cw:bili}.

\setcounter{equation}{0}

\section{Graded algebras with a single homogeneous relation}
Throughout this paper, $k$ is a field and $V$ is a $k$-vector space. The vector space $V$ is assumed to be graded by $V=V_1$, 
i.e., all elements of $V$ are of degree 1. The tensor algebra $T(V)$ of $V$ is graded accordingly. If $R$ 
is a graded subspace of $T(V)$, then the graded algebra $A=A(V,R)$ is defined as being the quotient of 
the tensor algebra $T(V)$ by the 2-sided ideal $I(R)$ generated by $R$. If $R$ is a subspace of $V^{\otimes N}$
for a certain integer $N\geq 2$, then $R$ is graded and the graded algebra $A(V,R)$ is said to be 
\emph{N-homogeneous}. 

If $M$ is a $k$-vector space, the set $\mathcal{L}(M)$ of vector subspaces of $M$, ordered by inclusion of sets, is 
a \emph{lattice}. Infimum and supremum in this lattice are respectively intersection and sum of subspaces. In general, a triple 
of subspaces of $M$ is not distributive, e.g., three distinct one-dimensional subspaces of a two-dimensional space. A sublattice 
of $\mathcal{L}(M)$ in which every triple of subspaces is distributive, is said to be \emph{distributive}. We extend this 
terminology to algebras as follows.
\Bdf
Let $A=A(V,R)$ be an $N$-homogeneous graded algebra. We say that $A$ is distributive if for any $m\geq 0$ ($m\geq N+2$
suffices), the sublattice of $\mathcal{L}(V^{\otimes m})$ generated by the subspaces $V^{\otimes
i}\otimes R\otimes V^{\otimes j}$ such that $i+N+j=m$, is distributive.
\Edf

Let $A=A(V,R)$ be an $N$-homogeneous graded algebra. Define the map $\nu _N : \mathbb{N}\rightarrow \mathbb{N}$ by 
$$\nu _N (2i)=Ni,$$
$$\nu _N (2i+1)= Ni+1,$$
for any $i \in \mathbb{N}$. For graded algebras where other functions of the degrees of generation of the projectives in the minimal resolution of $k$ are discussed, the reader is referred to~\cite{bbk:periodic, gm:delta, gmmz:dkos, gs:fingen}. 
 
For $m\geq 0$, introduce the subspace $W_{m}$ of $V^{\otimes m}$ by
$$W_{m}=\bigcap_{i+N+j=m}V^{\otimes i}\otimes R\otimes V^{\otimes j}.$$
Then $W_{m}=V^{\otimes m}$ if $m<N$, and $W_{N}=R$. As defined in~\cite{rb:nonquad}, the \emph{Koszul complex} $K(A)$
of the $N$-homogeneous graded algebra $A$ is the following complex 
\begin{equation}
\cdots \longrightarrow K_{i} \stackrel{\delta_{i}}{\longrightarrow} K_{i-1} \longrightarrow \cdots
\longrightarrow K_{1} \stackrel{\delta_{1}}{\longrightarrow} K_{0} \longrightarrow 0\,,
\end{equation}
in which $K_{i}=A\otimes W_{\nu _N(i)}$, and the left $A$-linear differential $\delta_{i}$ is defined by the inclusion
of $W_{\nu _N(i)}$ into $V^{\otimes (\nu _N(i)-\nu _N(i-1))}\otimes W_{\nu _N(i-1)}$. 
The homology of $K(A)$ is $k$ in degree 0, and 0 in degree 1.
\Bdf \label{defK}
An $N$-homogeneous graded algebra $A$ is said to be $N$-Koszul if the homology of $K(A)$ is $0$ in every degree $i>0$.
\Edf

For equivalent definitions, see~\cite{rb:nonquad}. The next proposition is essential for the sequel.
\Bpo \label{ess}
Assume that the $N$-homogeneous graded algebra $A$ is distributive. Then the following are equivalent:

\emph{(i)} $A$ is $N$-Koszul.

\emph{(ii)} $(R\otimes V^{\otimes m}) \cap (V^{\otimes m}\otimes R) \subseteq V^{\otimes (m-1)}\otimes R \otimes
V \ \ \mathrm{for} \ m=2, \ldots , N-1.$

\emph{(iii)} $(R\otimes V^{\otimes m}) \cap (V^{\otimes m}\otimes R) = W_{N+m} \ \ \mathrm{for} \ m=2, \ldots , N-1.$
\Epo
\Bdm
Proposition 3.4 and Lemma 2.6 in~\cite{rb:nonquad} provide respectively the equivalences
(i)$\Leftrightarrow$(ii) and (ii)$\Leftrightarrow$(iii).
\qed
\\ 

We are now ready to limit ourselves to graded algebras $A$ having a single homogeneous relation. If $N$ is the degree 
of this relation, it is clear that $A$ is $N$-homogeneous. Let us recall Gerasimov's theorem~\cite{gera:distrib}. As we said 
in the Introduction, it holds if $V$ is any positively graded vector space, i.e., 
$V=V_1\oplus V_2 \oplus \cdots $.
But assuming $V=V_1$ is natural as far as $N$-Koszul algebras are concerned, and it is our framework here. 
The question of defining $N$-Koszul algebras when generators are not necessarily of degree 1 deserves attention, 
and Gerasimov's theorem could be a guide for doing that.
\Bte \label{gera}
Let $A=A(V,R)$ be an $N$-homogeneous graded algebra. If $R$ is one-dimensional, then $A$ is distributive.
\Ete

Joining the two previous results, we get immediately the following criterion.
\Bte \label{crit2}
Let $A=A(V,R)$ be an $N$-homogeneous graded algebra. If $R$ is one-dimensional, then the following are equivalent.

\emph{(i)} $A$ is $N$-Koszul.

\emph{(ii)} $(R\otimes V^{\otimes m}) \cap (V^{\otimes m}\otimes R) \subseteq V^{\otimes (m-1)}\otimes R \otimes
V \ \ \mathrm{for} \ m=2, \ldots , N-1.$

\emph{(iii)} $(R\otimes V^{\otimes m}) \cap (V^{\otimes m}\otimes R) = W_{N+m} \ \ \mathrm{for} \ m=2, \ldots , N-1.$
\Ete

\setcounter{equation}{0}

\section{Antisymmetric relations}
A homogeneous element $f\in T(V)$ of degree $N$ is said to be \emph{antisymmetric} if 
$\sigma \cdot f=\mathrm{sgn} (\sigma) f$ for any $\sigma$ in the permutation group $S_N$, where $\sigma \cdot f$ 
denotes the natural action of $S_N$ on $V^{\otimes N}$ and sgn the sign. 
We also need a notation for the antisymmetrisers: if 
$v_1, \ldots , v_m$ are in $V$, we set
$$\mathrm{Ant}(v_1, \ldots , v_m)= \sum _{\sigma \in S_m} \mathrm{sgn}(\sigma)\, v_{\sigma (1)} \ldots v_{\sigma (m)}.$$
\Bpo \label{anti2}
Let $A=A(V,R)$ be an $N$-homogeneous graded algebra. Assume that $k$ has characterictic zero, $V$ is 
finite-dimensional of dimension $n\geq 2$, and $R$ is generated by an antisymmetric tensor $f\neq 0$ of degree $N$ 
with $2\leq N \leq n$. Then $A$ is $N$-Koszul and the global dimension of $A$ is $2$.
\Epo
\Bdm
The result is a consequence of the relations
\begin{equation} \label{equal}
(R\otimes V^{\otimes m}) \cap (V^{\otimes m}\otimes R) =0, \ \mathrm{for} \ m=1, \ldots , N-1.
\end{equation} 
Assuming these relations, $A$ is $N$-Koszul by Theorem \ref{crit2}. Furthermore (\ref{equal}) for $m=1$ implies that 
$W_m=0$ for every 
$m\geq N+1$, while $W_N\neq 0$. Thus the global dimension of $A$ is 2 (Corollary 2.6 in~\cite{rbnm:kogo}). It remains to 
prove (\ref{equal}). Fix $1\leq m \leq N-1$. An element of $(R\otimes V^{\otimes m}) \cap (V^{\otimes m}\otimes R)$ 
is antisymmetric because the $R$'s overlap. Then (\ref{equal}) is an immediate consequence of the following lemma. \qed
\Blm
Assume $1\leq s < r \leq n$. A nonzero antisymmetric tensor $a$ of degree $r$ is never divisible on the left (or on the right) 
by an antisymmetric tensor of degree $s$.
\Elm
\Bdm
Let $(x_1, \ldots, x_n)$ be a basis of $V$. Write uniquely
$$a= \sum _{1\leq i_1< \cdots < i_r \leq n} \alpha_{i_1\ldots i_r}\, \mathrm{Ant}(x_{i_1},\ldots , x_{i_r}).$$ 
There exist $i_1, \ldots ,i_r$ such that $\alpha_{i_1\ldots i_r}\neq 0$. Assume that
$$a = \left(\sum _{1\leq j_1< \cdots < j_s \leq n} \beta_{j_1\ldots j_s} \, \mathrm{Ant}(x_{j_1},\ldots , x_{j_{s}})\right)
\left(\sum _{1\leq k_1, \ldots , k_{r-s} \leq n} \gamma_{k_1\ldots k_{r-s}} \, x_{k_1} \ldots  x_{k_{r-s}} \right).$$
Examining on the RHS the coefficients of $x_{i_1}\ldots x_{i_r}$ and of 
$$x_{i_1}\ldots x_{i_{s-1}} x_{i_{s+1}} x_{i_{s}} x_{i_{s+2}} \ldots x_{i_r},$$
we deduce that 
$$\alpha_{i_1\ldots i_r}= \beta_{i_1\ldots i_s}\, \gamma_{i_{s+1}\ldots i_{r}} = - \beta_{i_1\ldots i_{s-1}i_{s+1}}\,
\gamma_{i_s i_{s+2}\ldots i_{r}},$$
thus $\beta_{i_1\ldots i_s}\neq 0$ and $\gamma_{i_s i_{s+2}\ldots i_{r}}\neq 0$. But the coefficient of 
$x_{i_1}\ldots x_{i_{s}} x_{i_{s}} x_{i_{s+2}} \ldots x_{i_r}$ is zero on the LHS, and it is equal to 
$\beta_{i_1\ldots i_s} \gamma_{i_s i_{s+2}\ldots i_{r}}$ on the RHS, hence we get a contradiction. 
\qed
\\

The Hilbert series of a graded algebra $A$ is denoted by $H_A(t)$. 
For any $N$-Koszul algebra $A$ such that the $A$-modules forming the Koszul complex of $A$ are of finite type (i.e., such that $V$ is finite-dimensional), it is known that 
\begin{equation} \label{hilb}
H_A(t) \left(\sum_{i\geq 0} (-1)^i \dim (W_{\nu_N(i)})\, t^{\nu_N(i)}\right)= 1.
\end{equation}
This formula is just a consequence of the Euler-Poincar\'e identity applied to the Koszul complex of $A$ 
(see e.g.~\cite{atv:modules, dvp:plac, kriegk:cras}), and it can be viewed as a duality formula between $A$ and the Yoneda algebra of $A$.
In some cases, the computation of $\dim (W_{\nu_N(i)})$ is possible for every $i$, and (\ref{hilb}) allows us 
to compute $\dim (A_i)$ for every $i$. It is particularly interesting when a linear basis of $A_i$ is known but 
counting the number $a_i$ of elements of this basis is tedious or unclear. This method to compute $a_i=\dim (A_i)$ 
is called the \emph{numerical Hilbert series combinatorics}. With the notations and assumptions of Proposition \ref{anti2}, 
it is easy to see that in this case $a_i$ is the number of monomials in $x_1, \ldots, x_n$ 
not containing $x_1 x_2\ldots x_N$ as factor. Since 
$H_A(t)= (1-nt + t^N)^{-1}$ by (\ref{hilb}), we get directly the recursive formula $a_i= n a_{i-1} - a_{i-N}$ with 
$a_0=1$ and $a_i=0$ when 
$i<0$. Moreover, as all the complex zeros of the polynomial $1-nt + t^N$ 
have module 1 if and only if $n=2$, the Gelfand-Kirillov dimension of 
$A$ is 2 if $n=2$, and $\infty$ if $n>2$ (see Proposition 2.13 in~\cite{atv:modules}). 

Recall the definition of AS-Gorenstein algebras and the definition of Calabi-Yau algebras. Here AS stands for Artin and Schelter~\cite{as:regular}. The definition of Calabi-Yau algebras is due to Victor Ginzburg~\cite{g:cy}. In these definitions, the functor Ext is applied to left modules.
\Bdf
(i) Let $A$ be a connected graded $k$-algebra of finite global dimension $n\geq 1$. We say that $A$ is
AS-Gorenstein if \emph{Ext}$_{A}^{i}(k,A)=0$ for $i \neq n$, and \emph{Ext}$_{A}^{n}(k,A)\cong k$ as right $A$-modules. 

(ii) Let $A$ be a $k$-algebra having a finite projective $A$-bimodule resolution by bimodules of finite type. We say that $A$ is a Calabi-Yau algebra of dimension $n\geq 1$ (or a $n$-Calabi-Yau algebra) if \emph{Ext}$_{A^e}^{i}(A,A^e)=0$ for $i \neq n$, and \emph{Ext}$_{A^e}^{n}(A,A^e)\cong A$ as right $A^e$-modules, where $A^e=A\otimes A^{op}$ and $A^{op}$ is the opposite algebra of $A$.
\Edf
It is often possible in the case of an $N$-Koszul algebra $A$ of finite global dimension to decide whether $A$ is AS-Gorenstein or 
not (\cite{rbnm:kogo}, Theorem 5.4). Keeping again the notations and 
assumptions of Proposition \ref{anti2}, it is clear from Proposition 5.2 in~\cite{rbnm:kogo} that if $A$ is AS-Gorenstein, then 
$N=2$. Assume that $N=2$. Let $(x_1, \ldots, x_n)$ be a basis of $V$. Write 
$f=\sum_{1\leq i<j\leq n} \lambda_{ij}\, \mathrm{Ant}(x_i,x_j)$ and denote by $M$ the antisymmetric matrix with coefficients $\lambda_{ij}$. 
According to Theorem 3 in~\cite{mdv:multi}, $A$ is AS-Gorenstein if and only if $M$ is invertible (see also Section 5). 
If $M$ is invertible, then $n$ is even and $f$ is the standard symplectic form in an appropriate basis, so that $A$ is the preprojective algebra of a non-Dynkin quiver, hence $A$ is 2-Calabi-Yau~\cite{bock:gcy} (see 7.1 in~\cite{boc:nccone} for an interpretation of $A$ as a noncommutative tangent cone). Summing up now what we have obtained.
\Bpo
Let us keep the notations and assumptions of \emph{Proposition \ref{anti2}}. Then $A$ is AS-Gorenstein if and only if $N=2$, $n$ is even 
and $f$ (viewed as an antisymmetric bilinear form) is nondegenerate. When $N=2$ and $n$ is even, 
there is exactly one AS-Gorenstein algebra $A$ up to a linear change of generators, and this algebra $A$ is $2$-Calabi-Yau.
\Epo

If we have $N=n$ in the statement of Proposition \ref{anti2}, $R$ coincides with the subspace of all antisymmetric 
tensors of degree $N$. It was already known that $A(V,R)$ is $N$-Koszul in this case, since $A(V,R)$ is $N$-Koszul 
if $R$ is the subspace of all antisymmetric tensors of degree $N$ for any $2\leq N \leq n$~\cite{rb:nonquad}. 
This suggests to deduce the $N$-Koszulity of the algebra $A(V,R)$ of Proposition \ref{anti2} when $2\leq N < n$ 
from the known case $N=n$, 
by adding new variables. Actually, this method is possible only if the antisymmetric tensor generating $R$ is an antisymmetriser. In order to apply this method, we need some considerations on the free product of algebras.

Let $A=A(V,R)$ and $A'=A(V',R')$ be $N$-homogeneous graded algebras. Their \emph{free product} 
is the $N$-homogeneous algebra $A\ast A'$ whose space of generators is $V\oplus V'$ 
and space of relations is $R\oplus R'$~\cite{bafr:kalg}. Here $R \oplus R'$ is considered in
$(V\oplus V')^{\otimes N}$ via the canonical injective map from $V^{\otimes N} \oplus  V^{'\otimes N}$ 
into $(V\oplus V')^{\otimes N}$. If $A$ and $A'$ are distributive, 
then $A\ast A'$ is distributive (Theorem 1 (d) in~\cite{bafr:kalg}). 
\Bpo
Let $A$ and $A'$ be $N$-homogeneous and distributive. If $A$ and $A'$ are $N$-Koszul, then
$A\ast A'$ is $N$-Koszul.
\Epo
\Bdm
Set $\mathcal{V}=V\oplus V'$. Define $\mathcal{R}$ as being $R\oplus R'$ naturally viewed in
$\mathcal{V}^{\otimes N}$. Fix $2\leq m \leq N-1$. 
Omit tensor symbols $\otimes$ in order to simplify notations. Let $a$ belong to 
$(\mathcal{R}\mathcal{V}^m) \cap (\mathcal{V}^m \mathcal{R})$. 
Using Proposition \ref{ess}, it suffices to show that $a$ belongs to $\mathcal{V}^{m-1} \mathcal{R} \mathcal{V}$. 
Write 
$$a=a_1+a_2=b_1+b_2,$$
where $a_1 \in R\mathcal{V}^m$, $a_2\in R'\mathcal{V}^m$, $b_1\in \mathcal{V}^m R$ and 
$b_2\in \mathcal{V}^m R'$. Set $\mathcal{V}_1=V$ and  $\mathcal{V}_2=V'$. The vector space 
$\mathcal{V}^{N+m}$ is the direct sum of the subspaces $\mathcal{V}_{i_1}\ldots \mathcal{V}_{i_{N+m}}$ for all 
$(N+m)$-tuples $(i_1, \ldots ,i_{N+m})$ of elements of the set $\{1,2\}$. Using this direct sum and the fact that 
$a_1 \in V^N\mathcal{V}^m$, $a_2\in V'^N\mathcal{V}^m$, $b_1\in \mathcal{V}^m V^N$ and  
$b_2\in \mathcal{V}^m V'^N$, we see that $b_1$ is in $V^m V^N$ or in $V'^m V^N$. 
But if $b_1 \in V'^m V^N$, then $b_1=a_2$ which implies that $b_1=a_2=0$ (since $m<N$), 
hence $a=a_1=b_2$ which implies again $a=0$. Assuming that $a$ is nonzero, we have shown that 
$b_1$ is in $V^m V^N$ and that $b_1=a_1$. Thus $a_1$ belongs to $(RV^m) \cap (V^m R)$. 
For the same reasons, $a_2$ belongs to $(R'V'^m) \cap (V'^m R')$. As $A$ and $A'$ are Koszul, 
$a_1$ belongs to $V^{m-1}RV$ and $a_2$ belongs to $V'^{m-1}R'V'$. Thus $a=a_1+a_2$ belongs to 
$\mathcal{V}^{m-1} \mathcal{R} \mathcal{V}$.
\qed
\\

Assume $k$ is of characteristic zero, $V$ is of finite dimension $n\geq 3$ and $(x_1, \ldots, x_n)$ is a basis of $V$. 
Fix $2\leq N <n$. Let $A$ be the $N$-homogeneous graded algebra defined by the generators $x_1, \ldots, x_n$ and 
the single relation $\mathrm{Ant}(x_1,\ldots , x_N)=0$. Let $B$ be defined by the generators 
$x_1,\ldots , x_N$ only and the same relation. Then $A$ is the free product of $B$ and $B'$, where $B'$ is the free algebra in the generators 
$x_{N+1}, \ldots, x_n$. We know that $B$ is distributive and $N$-Koszul~\cite{rb:nonquad} (but this fact is not trivial), 
and it is clear for $B'$. Thus we conclude that $A$ is distributive and $N$-Koszul, without using Gerasimov's theorem. Let us notice that this method is not general since an antisymmetric tensor $f$ of degree $N$, $2\leq N <n$, is not in general an antisymmetriser even up to a linear change of generators, for example if $f$ is a symplectic form over $n\geq 4$ variables. 

\setcounter{equation}{0}

\section{Monomial relations}
Throughout this section, $V$ denotes a $k$-vector space of finite dimension $n\geq 1$, and $(x_1, \ldots, x_n)$ a basis of $V$. 
We need some special terminology. An element $a$ of $T(V)$ of the form $a=x_{i_1}\ldots x_{i_r}$ for some $r\geq 1$ is 
called a \emph{monomial of degree} $r$. If $1\leq s \leq r$, the monomials of the form $x_{i_u}\ldots x_{i_{u+s-1}}$ 
with $1\leq u \leq r-s+1$ are called the \emph{factors of degree s} of $a$. Moreover $u$ is called the \emph{index} of the 
factor $x_{i_u}\ldots x_{i_{u+s-1}}$ of degree $s$ of $a$, and the factor of index 1 (resp. $r-s+1$) is called the 
first (last) factor of degree $s$ of $a$. The next proposition is identical to Proposition 3.8 in~\cite{rb:nonquad}. 
\Bpo \label{mono}
Fix an integer $N\geq 2$. Let $\mathcal{C}$ be a set of monomials of degree $N$. Let $R$ be the subspace of 
$V^{\otimes N}$ generated by $\mathcal{C}$, and $A=A(V,R)$ the $N$-homogeneous graded algebra associated to $V$ and $R$ (such an algebra is said to be \emph{monomial}). 
The following are equivalent.

\emph{(i)} $A$ is $N$-Koszul.

\emph{(ii)} If a monomial $a$ of degree $N+m$ where $2\leq m \leq N-1$ is such that 
its first factor and its last factor of degree $N$ belong to $\mathcal{C}$, then its factor of degree $N$ and 
index $m$ belongs to $\mathcal{C}$.

\emph{(iii)} If a monomial $a$ of degree $N+m$ where $2\leq m \leq N-1$ is such that its first factor and its last factor of 
degree $N$ belong to $\mathcal{C}$, then all its factors of degree $N$ belong to $\mathcal{C}$.
\Epo

A monomial $N$-homogeneous graded algebra $A$ is distributive, since distributivity of $A$ is reduced here to 
distributivity of $\cap$ over $\cup$. Thus Proposition \ref{mono} is an immediate consequence of Proposition \ref{ess}. 
Proposition \ref{mono} shows that any monomial 2-homogeneous graded algebra is Koszul, 
recovering a result by Fr\"oberg~\cite{fro:mono}. 

We are now interested in monomial algebras defined by a single monomial relation. 
In~\cite{rb:nonquad}, the following example is given: if $N=3$ and $\mathcal{C}=\{a\}$, then $A$ is 3-Koszul 
if and only if $a$ is not of the form $x_ix_jx_i$ with $i\neq j$. More generally, for any value of $N$, 
the next proposition provides 
the singletons $\mathcal{C}$ which are $N$-Koszul (i.e., such that $A$ is $N$-Koszul).
\Bpo \label{onemon}
Let $V$ be a $k$-vector space of finite dimension $n\geq 1$. Fix a basis $(x_1, \ldots, x_n)$ of $V$. 
Let $f=x_{i_1} \ldots x_{i_N}$ be a monomial of degree $N\geq 2$, and $R$ be the subspace of $V^{\otimes N}$ 
generated by $f$. Then $A=A(V,R)$ is $N$-Koszul if and only if there exists no $m$ 
in $\{2,\ldots , N-1\}$ such that 
$$f=(x_{i_1} \ldots x_{i_m})^q\, x_{i_1} \ldots x_{i_r},$$
where $N=mq+r$ is the division of $N$ by $m$ with remainder $r$, and where $i_1,\ldots ,i_m$ are not all equal.
\Epo
Before proving Proposition 4.2, we need to prove the following lemma.
\Blm \label{lemw}
Keep the notations and assumptions of \emph{Proposition \ref{onemon}}. If $f=x_i^N$ for some $i$, then $W_m=k\cdot x_i^m$ for every $m\geq N+1$. 
Otherwise, $W_m=0$ for every $m\geq N+1$.
\Elm
\Bdm
Let $a$ be nonzero in $W_{N+1}$. Write 
$$a=x_{i_1} \ldots x_{i_N} \left(\sum_{1\leq i \leq N} \lambda_i x_i \right) = 
\left(\sum_{1\leq j \leq N} \mu_j x_j \right) x_{i_1} \ldots x_{i_N},$$
where $\lambda_i \neq 0$ for a certain $i$. Then there exists $j$ such that $x_jx_{i_1} \ldots x_{i_N}=
x_{i_1} \ldots x_{i_N}x_i$, hence $j=i_1,\, i_1=i_2, \ldots ,i_{n-1}=i_N,\, i_N=i$. Thus $f=x_i^N$. 
Therefore, if $f$ is not of type $x_i^N$, then $W_{N+1}=0$, hence $W_m=0$ for every $m\geq N+1$. 
If $f=x_i^N$, it is obvious that $W_m=k\cdot x_i^m$ for every $m\geq N+1$.
\qed
\\

We are now ready to prove Proposition \ref{onemon}. If $f=x_i^N$ for some $i$, then Proposition \ref{mono} 
shows that $A$ is $N$-Koszul. For the sequel, assume that $f$ is not of the form $x_i^N$ for some $i$. 
Lemma \ref{lemw} and Theorem \ref{crit2} (iii)
show that $A$ is $N$-Koszul if and only if for every $m\in \{2, \ldots ,N-1\}$ there is no overlap monomial for $f$ of degree $N+m$, 
meaning that there is no monomial of degree $N+m$ whose first and last factor of degree $N$ is equal to $f$. 

Let $m\in \{2, \ldots ,N-1\}$ such that $f=(x_{i_1} \ldots x_{i_m})^q\, x_{i_1} \ldots x_{i_r}$, where 
$N=mq+r$ with $0 \leq r \leq N-1$, and $i_1,\ldots ,i_m$ are not all equal (hence $f\neq (x_{i_1})^ N$). Then
$$f x_{i_{r+1}} \ldots x_{i_m} x_{i_1} \ldots x_{i_r} = x_{i_1} \ldots x_{i_m} f$$
is an overlap monomial for $f$ of degree $N+m$. Thus $A$ is not $N$-Koszul.

Conversely, assume that $A$ is not $N$-Koszul. Let $m$ be in $\{2, \ldots ,N-1\}$ such that there exist $j_1,\ldots , j_m$ 
and $k_1,\ldots k_m$ in $\{1, \ldots ,n\}$ satisfying
$$x_{i_1} \ldots x_{i_N} x_{j_1} \ldots x_{j_m}=x_{k_1} \ldots x_{k_m} x_{i_1} \ldots x_{i_N}.$$
Then $i_1=k_1,\ldots, i_m=k_m, i_{m+1}=i_1, \ldots , i_{2m}=i_m, \ldots, i_{mq}=i_{m(q-1)}, \ldots, i_N=i_{m(q-1)+r}, 
j_1=i_{m(q-1)+r+1}, \ldots , j_m=i_{mq+r}=i_N$, thus $f=(x_{i_1} \ldots x_{i_m})^q\, x_{i_1} \ldots x_{i_r}$. 
Since $f\neq (x_{i_1})^ N$, $i_1,\ldots ,i_m$ are not all equal. 
\qed
\\

The proof of the following is left to the reader.
\Bpo \label{onemon1}
Let us keep the notations and assumptions of \emph{Proposition \ref{onemon}}. Assume that $A$ is $N$-Koszul. 
If $f=x_i^ N$ for some $i$, then the global dimension of $A$ is infinite, we have
$$H_A(t)=(1-nt+t^ N-t^{N+1}+t^{2N}-t^ {2N+1}+ \cdots )^ {-1},$$
and the Gelfand-Kirillov dimension of $A$ is $0$ if $n=1$, $\infty$ if $n>1$. 
Otherwise, the global dimension of $A$ is equal to $2$, we have
$$H_A(t)=(1-nt+t^ N)^ {-1},$$
and the Gelfand-Kirillov dimension of $A$ is $2$ if $n=N=2$, $\infty$ if $n>2$ or $N>2$.
In either case, $A$ is never AS-Gorenstein. 
\Epo
\addtocounter{Df}{1}
\Rm 
There is no available description (analogous to that of Proposition \ref{onemon}) of the $N$-Koszul sets $\mathcal{C}$ when $n \geq 2$ and the number $p$ of elements of $\mathcal{C}$ is fixed with $1<p<n^N$. The number of such ($N$-Koszul or not $N$-Koszul) $\mathcal{C}$ is equal to the binomial coefficient $\binom{n^N}{p}$, and it is not clear how to find an algorithm which is polynomial in $N$ and which lists all the $N$-Koszul $\mathcal{C}$'s when $n \geq 2$ and $1<p<n^N$. Basing ourselves on some computer evidence, we conjecture that the number Koszul$(n,N,p)$ of these $N$-Koszul $\mathcal{C}$'s is zero when $n^N\! /2 <p < n^N$ and $N$ is large enough.

\setcounter{equation}{0}

\section{Gerasimov's Theorem for $N=2$}
Let us start with the following well-known result (see Proposition 10.1 in~\cite{gmmz:dkos} for more general assumptions including quiver algebras with relations).
\Bpo \label{gld2}
Let $A=A(V,R)$ be an $N$-homogeneous graded algebra. If the global dimension of $A$ is $2$, then $A$ is $N$-Koszul.
\Epo
\Bdm
Clearly $V$ is a minimal space of generators and $R$ is a minimal space of relations. Then the beginning 
\begin{equation} \label{debres1}
A\otimes R \stackrel{\delta_{2}}{\longrightarrow} A\otimes V \stackrel{\delta_{1}}{\longrightarrow} A \longrightarrow 0
\end{equation}
of the Koszul complex $K(A)$ of $A$, with the natural projection $\epsilon : A\rightarrow k$, is the beginning of a minimal projective resolution of $k$ in the category $A$-grMod (whose objects are the graded left $A$-modules and morphisms are the homogeneous $A$-linear maps of degree 0). Since the global dimension of $A$ is equal to the length of a minimal projective resolution of $k$, then 
\begin{equation} \label{res1}
0 \longrightarrow A\otimes R \stackrel{\delta_{2}}{\longrightarrow} A\otimes V \stackrel{\delta_{1}}{\longrightarrow} A \longrightarrow 0
\end{equation}
is a minimal projective resolution of $k$. But $W_{N+1}$ is contained in $\ker(\delta_{2})$, thus vanishes, and therefore (\ref{res1}) coincides with $K(A)$. 
\qed
\\

We recall another basic result due to Backelin~\cite{bafr:kalg}. For a proof, see Section 3 in~\cite{rb:nonquad} or Theorem 4.1 in~\cite{pp:quad}. The fact that distributivity of $A$ implies Koszulity of $A$ is a consequence of Proposition 2.3 with $N=2$. The converse is the part which is hard to prove.
\Bpo \label{back}
Let $A=A(V,R)$ be an $2$-homogeneous graded algebra. Then $A$ is Koszul if and only if $A$ is distributive.
\Epo

From now on and throughout this section, $V$ denotes a $k$-vector space of finite dimension $n\geq 1$. Let $f\neq 0$ be in $V\otimes V$ and $R$ be the one-dimensional subspace of $V\otimes V$ generated by $f$. We are interested in the 2-homogeneous graded algebra $A=A(V,R)$. If $(x_1, \ldots, x_n)$ is a basis of $V$, the $n\times n$ matrix $M=M(f)=(f_{ij})_{1\leq i,j\leq n}$ where
$$f=\sum_{1\leq i,j\leq n} f_{ij}\, x_i x_j,$$
is called the \emph{matrix of coefficients} of $f$ relative to this basis.
Denote by $x$ the $n\times 1$ matrix with entries $x_1, \ldots, x_n$ and $x^t=(x_1, \ldots, x_n)$ its transpose. Denote by $A^{1\times n}$ (resp. $A^{n\times 1}$) the left (resp. right) $A$-module of the $1\times n$ (resp. $n\times 1$) matrices with entries in $A$.  There is an obvious isomorphism from (\ref{debres1}) to
\begin{equation} \label{debres2}
A \stackrel{\cdot (x^tM)}{\longrightarrow} A^{1\times n} \stackrel{\cdot x}{\longrightarrow} A \longrightarrow 0
\end{equation}
which is therefore the beginning of a minimal projective resolution of $k$ in $A$-grMod with the same arrow $\epsilon : A\rightarrow k$.
\Blm \label{lem}
Assume that the global dimension of $A$ is $2$. Then $A$ is Koszul. Moreover, $A$ is AS-Gorenstein if and only if $M$ is invertible.
\Elm
\Bdm
The first assertion comes from Proposition \ref{gld2}. Notice that $n\geq 2$ since $n=1$ implies that the global dimension of $A$ is infinite. Applying the functor Hom$_A(-,A)$ to the minimal projective resolution
\begin{equation} \label{res2}
0 \longrightarrow A \stackrel{\cdot (x^tM)}{\longrightarrow} A^{1\times n} \stackrel{\cdot x}{\longrightarrow} A \longrightarrow 0,
\end{equation}
we get 
\begin{equation} \label{res3}
0 \longrightarrow A \stackrel{x \cdot}{\longrightarrow} A^{n\times 1} \stackrel{(x^tM) \cdot}{\longrightarrow} A \longrightarrow 0.
\end{equation}
Saying that $A$ is AS-Gorenstein means that (\ref{res3}) is a resolution of $k$ in grMod-$A$ (whose objects are the graded right $A$-modules and morphisms are the homogeneous $A$-linear maps of degree 0) for a certain arrow $A\rightarrow k$ of grMod-$A$. Actually, as the latter arrow is necessarily of the form $\lambda \epsilon$ with $\lambda \neq 0$ in $k$, we can impose that this arrow is $\epsilon$. Consider now the following commutative diagram in the category grMod-$A$:
$$\xymatrix@=1cm{0\ar[d]^{id}\ar[r] &A \ar[d]^{id}\ar[r]^{x\cdot}&A^{n\times 1}\ar[d]^{M\cdot}\ar[r]^{(x^tM) \cdot}&A\ar[d]^{id}\ar[r]^{\epsilon}& k \ar[d]^{id}\ar[r] & 0 \ar[d]^{id} \\
0\ar[r]&A\ar[r]^{(Mx)\cdot}&A^{n\times 1}\ar[r]^{x^t \cdot}&A\ar[r]^{\epsilon}& k \ar[r] & 0  }$$

The second row is exact because it is isomorphic to the Koszul resolution of $k$ in grMod-$A$, and it is a basic fact that Definition \ref{defK} is equivalent to its right version. If $M$ is invertible, then the first row is exact, thus $A$ is AS-Gorenstein. Assume that $A$ is AS-Gorenstein. Then the first row is a projective resolution of $k$ in grMod-$A$, and the point is that this resolution is \emph{minimal} since from right to left: $A$ is generated in degree 0, $A^{n\times 1}$ is generated in degree 1, $A$ is generated in degree 2, these degrees are strictly increasing, so that when applying the functor $-\otimes_A k$ the morphisms of the first row vanish. It is again a basic fact of homological algebra that a morphism between two minimal projective resolutions is an isomorphism (\!~\cite{bou:alghom}, 3, no. 6, Proposition 8), hence $M$ is invertible.
\qed
\\

The rank of the matrix of coefficients of $f$ relative to a basis of $V$ does not depend on the basis, and it is called the \emph{rank} of $f$. The rank of $f$ is $\geq 1$, and it is equal to 1  if and only if the following property (P1) holds:

\noindent
(P1) there exists a basis of $V$ for which $M(f)$ has only one nonzero column.

We also need the stronger property (P2) which is strictly stronger whenever $n\geq 2$:

\noindent
(P2) there exists a basis of $V$ for which all the entries of $M(f)$ are zero, except one which is located on the diagonal.

Notice that (P2) is equivalent to say that $f$ is symmetric of rank 1. The following contains Theorem 3 in~\cite{mdv:multi}.
\Bpo \label{nP1}
Assume that the rank of $f$ is $>1$. Then

\emph{(i)} every nonzero element of $V$ is neither a left nor a right zerodivisor in $A$,

\emph{(ii)} the global dimension of $A$ is $2$,

\emph{(iii)} $A$ is Koszul,

\emph{(iv)} $A$ is AS-Gorenstein if and only if $M$ is invertible.
\Epo
\Bdm
To prove (i), it is sufficient to prove that for every nonzero $v\in V$ and every $a\in A_i$ ($i\geq 0$), $av=0$ implies that $a=0$. Proceed by induction on $i$, the case $i=0$ being clear. Assume that the property holds for some $i\geq 0$. Let $v$ be nonzero in $V$ and $a\in A_{i+1}$ be such that $av=0$. Choose a basis $(y_1, \ldots , y_n)$ of $V$ such that $y_1=v$, and write $f=\sum_{1\leq i,j\leq n} f'_{ij}\, y_i y_j$. Then $(a,0 \ldots , 0).y=0$ where $y$ denotes the column matrix of the $y_i$'s. Since the complex
\begin{equation} \label{res4}
0 \longrightarrow A \stackrel{\cdot (y^tM')}{\longrightarrow} A^{1\times n} \stackrel{\cdot y}{\longrightarrow} A \longrightarrow 0,
\end{equation}
is exact at $A^{1\times n}$ (where $M'=(f'_{ij})_{1\leq i,j\leq n}$), there exists $b\in A_{i-1}$ such that $b.\sum_{1\leq i\leq n} f'_{i1} y_i=a$ and $b.\sum_{1\leq i\leq n} f'_{ij} y_i=0$ if $j\neq 1$. The second condition is not void because $n\geq 2$. Since (P1) does not hold, there exists $j\neq 1$ such that $\sum_{1\leq i\leq n} f'_{ij} y_i\neq 0$. The induction hypothesis implies thus $b=0$, hence $a=0$.

From (i), we deduce that (\ref{res4}) is exact. In fact, if $a\in A$ is such that $a\cdot (y^tM')=0$, then $a.\sum_{1\leq i\leq n} f'_{ij} y_i=0$  for $j=1, \ldots , n$. As $M\neq 0$, there exists $j$ such that $\sum_{1\leq i\leq n} f'_{ij} y_i\neq 0$, and thus $a=0$ by (i). So (\ref{res4}) is a minimal projective resolution of $k$, hence (ii) holds and (iii), (iv) follow from Lemma \ref{lem}.
\qed

\Bpo \label{P1}
Assume that the rank of $f$ is equal to $1$. Then

\emph{(i)} $A$ is Koszul,

\emph{(ii)} the global dimension of $A$ is $2$ or infinite,

\emph{(iii)} the global dimension is infinite if and only if \emph{(P2)} holds,

\emph{(iv)} $A$ is not AS-Gorenstein.
\Epo  
\Bdm
There exists a basis $(y_1, \ldots , y_n)$ of $V$ such that $f=(\sum_{1\leq i\leq n} f'_{i1} y_i)y_1$. Denote by $Y$ the ordered basis $y_1< \cdots < y_n$ of $V$. We show that $A$ is $Y$-confluent, which will be enough to conclude that $A$ is Koszul~\cite{rb:conkos}. Let $i_0$ be the largest $i$ such that $f'_{i1}\neq 0$. The relation of $A$ is written
$$y_{i_0}y_1=-\frac{1}{f'_{i_0 1}} \sum_{1\leq i\leq n} f'_{i1} y_i y_1,$$
with $y_i y_1<y_{i_0} y_1$, so that $y_{i_0} y_1$ is the unique nonreduced monomial. Examine the ambiguous monomials of degree 3, i.e., the monomials $y_p y_q y_r$ with $y_p y_q$ and $y_q y_r$ nonreduced. If $i_0\neq 1$, there is no such $y_p y_q y_r$, hence $A$ is $Y$-confluent. If $i_0=1$, the single ambiguous monomial is $y_1^3$, the relation becomes $y_1^2=0$, and $A$ is still $Y$-confluent. Thus (i) is proved.

According to the Koszul resolution $K(A)$ of $A$, the global dimension of $A$ is 2 if and only if $(R\otimes V)\cap (V\otimes R) =0$. It is easy to check that (P2) holds if and only if $(R\otimes V)\cap (V\otimes R) \neq 0$. Therefore, if $A$ does not satisfy (P2), then the global dimension of $A$ is 2. If $A$ satisfies (P2), then the relation of $A$ is of the form $y_1^2=0$, so that the global dimension is infinite (Proposition \ref{onemon1}). Thus (ii) and (iii) are proved.

To prove (iv), we can assume that the global dimension is finite (by definition, an AS-Gorenstein algebra has a finite global dimension), hence is equal to 2, and we conclude by Lemma \ref{lem} since $n\geq 2$ and $M$ is not invertible.
\qed
\\

Under the hypothesis of Proposition \ref{nP1} or of Proposition \ref{P1}, we have $H_A(t)= (1-nt + t^2)^{-1}$ if the global dimension is 2, and $H_A(t)=(1-nt+t^2-t^3 + \cdots )^ {-1}$ if the global dimension is infinite. The Gelfand-Kirillov dimension of $A$ is derived as in Proposition \ref{onemon1} with $N=2$: it is equal to 0 if $n=1$, to $\infty$ if $n>2$, and if $n=2$, it is equal to the global dimension. Furthermore, joining Proposition \ref{nP1}, Proposition \ref{P1} and the hard part of Backelin's theorem (Proposition \ref{back}), we get Gerasimov's theorem for $N=2$.
\Bte
Let $V$ be a $k$-vector space of finite dimension $n\geq 1$. Let $f\neq 0$ be in $V\otimes V$ and $R$ be the one-dimensional subspace of $V\otimes V$ generated by $f$. Then the $2$-homogeneous graded algebra $A(V,R)$ is distributive.
\Ete
Recall that the $N$-homogeneous graded algebras $A(V,R)$ and $A(V',R')$ are said to be \emph{isomorphic} if there exists a bijective linear map $\phi : V \rightarrow V'$ such that $\phi ^{\otimes N} (R)=R'$ ~\cite{rbdvmw:homog}. We want to describe briefly the isomorphism classes of the $A(V,R)$'s when $V$ is a $\mathbb{C}$-vector space of given finite dimension $n\geq 1$ and $R$ is one-dimensional in $V\otimes V$. A first invariant for these isomorphism classes is the \emph{rank} (taking any integral value from 1 to $n$), i.e., the rank of $f$ generating $R$. The isomorphism class is said to be \emph{symmetric}, \emph{antisymmetric} or \emph{nondegenerate} if the matrix $M(f)$ is symmetric, antisymmetric or invertible. The isomorphism classes of AS-Gorenstein algebras are exactly the nondegenerate ones. There is only one isomorphism class containing the algebras having an infinite global dimension (respectively which are Calabi-Yau with necessarily $n$ even), and this class is symmetric of rank 1 (resp. antisymmetric of rank $n$). If $n\neq 2$, there is no isomorphism class formed of AS-regular algebras (i.e., AS-Gorenstein of finite Gelfand-Kirillov dimension). If $n=2$, the isomorphism classes of AS-regular algebras are parametrised by $\{0, 1\} \cup (\mathbb{C} \setminus \{0,1\})/(q \sim q^{-1})$, and there is only one symmetric (antisymmetric) such a class (Section 2 in~\cite{mdv:multi}). It is easy to see that all the isomorphism classes are parametrised by the set of orbits of all $n\times n$ complex matrices $M$ for the GL$(n,\mathbb{C})$-action of $P$ over $M$ given by $P^tMP$. The latter set has a long history, culminating recently to an explicit list of representatives of the orbits (see~\cite{cw:bili} and references therein).

\setcounter{equation}{0}

\section{Poincar\'e-Birkhoff-Witt deformations}
Let $A=A(V,R)$ be an $N$-Koszul graded algebra, $N\geq 2$. The tensor algebra $T(V)$ is filtered by the subspaces $F^n= k\oplus V \oplus \cdots \oplus V^{\otimes n}$, $n\in \mathbb{N}$. Let $\varphi : R \rightarrow F^{N-1}$ be $k$-linear. Define the subspace $P=(\mathrm{Id}-\varphi)(R)$ of $F^N$. Since $\mathrm{Id}-\varphi : R\rightarrow P$ is a linear isomorphism, $P$ and $R$ have the same dimension (which may be infinite). Then $U=T(V)/I(P)$ is a filtered algebra and there is a natural graded algebra morphism $p$ from $A$ onto the associated graded algebra gr$(U)$. The filtered algebra $U$ is said to be a \emph{Poincar\'e-Birkhoff-Witt (PBW) deformation} of $A$ if $p$ is an isomorphism~\cite{fv:pbwdef}.

Conversely, let $P$ be a subspace of $T(V)$. The filtered algebra $U=T(V)/I(P)$ is said to be Koszul if there exists an integer $N\geq 2$ such that $P\subseteq F^N$ and $U$ is a PBW deformation of $A=A(V,R)$ assumed to be Koszul, where $R=\pi (P)$ and $\pi : F^N \rightarrow V^{\otimes N}$ is the natural projection~\cite{bg:hsra}. Since $N$ is unique, we say that $U$ is $N$-Koszul. The first proposition determines the PBW deformations of $A$ when $A$ has a single monomial relation and is Koszul of infinite global dimension (see Proposition \ref{onemon1}).
\Bpo \label{infinite}
Assume that $V$ has finite dimension $n\geq 1$ and that $(x_1, \ldots, x_n)$ is a basis of $V$. Let $R$ be the subspace of $V^{\otimes N}$ generated by $x_1^N$, $N\geq 2$, and let $\varphi : R \rightarrow F^{N-1}$ be $k$-linear. Define the subspace $P=(\mathrm{Id}-\varphi)(R)$ of $F^N$. The filtered algebra $U=T(V)/I(P)$ is a PBW deformation of $A=A(V,R)$ if and only if $\varphi (x_1^N)$ is a polynomial in $x_1$.
\Epo
\Bdm
The space $W_{N+1}=(R\otimes V)\cap (V\otimes R)$ is generated by $x_1^{N+1}$. Write $\varphi = \sum_{j=0}^{N-1}\varphi_j$, $\varphi_j : R\rightarrow
V^{\otimes j}$. Set $\varphi_{j}^1= \varphi_j \otimes 1_V$ and $\varphi_{j}^2= 1_V \otimes \varphi_j$. According to the PBW theorem (\!\!~\cite{bg:hsra}, Theorem 3.4 and Proposition 3.6), $U$ is a PBW deformation of $A$ if and only if the following properties hold :
\begin{equation} \label{J'1}
(\varphi_{N-1}^1 - \varphi_{N-1}^2)(W_{N+1})\subseteq R,
\end{equation}
\begin{equation} \label{J'2}
\left(\varphi_j(\varphi^1_{N-1}-\varphi^2_{N-1})+\varphi^1_{j-1}-\varphi^2_{j-1}
\right)(W_{N+1})=0,\ 1\leq j \leq N-1,
\end{equation}
\begin{equation} \label{J'3}
\varphi_0(\varphi^1_{N-1}-\varphi^2_{N-1})(W_{N+1})=0.
\end{equation}
Clearly, (\ref{J'1}) holds if and only if $\varphi_{N-1}(x_1^N)=\lambda_{N-1} x_1^{N-1}$ for a certain $\lambda_{N-1}\in k$, and (\ref{J'1}) is equivalent to $(\varphi_{N-1}^1 - \varphi_{N-1}^2)(W_{N+1})=0$. If (\ref{J'1}) holds, then the conjunction of (\ref{J'2}) and (\ref{J'3}) is equivalent to saying that 
$$\left(\varphi^1_{j}-\varphi^2_{j} \right)(W_{N+1})=0,\ 0\leq j \leq N-2,$$
i.e., that $\varphi_j(x_1^N)=\lambda_{j} x_1^{j}$ for a certain $\lambda_{j}\in k$ for every $j=0, \ldots, N-2$.
\qed
\Bpo \label{gldim2}
Let $A=A(V,R)$ be an $N$-Koszul graded algebra, $N\geq 2$. Let $\varphi : R \rightarrow F^{N-1}$ be $k$-linear, and $U=T(V)/I(P)$ where $P=(\mathrm{Id}-\varphi)(R)$. If the global dimension of $A$ is $2$, then $U$ is a PBW deformation of $A$.
\Epo
\Bdm
The assumptions imply that $W_{N+1}=0$, thus (\ref{J'1})-(\ref{J'3}) hold.
\qed
\\

We want to illustrate Proposition \ref{gldim2} in the context of the previous section (assumed up to the end of the paper): $V$ has finite dimension $n\geq 1$, $f$ is a nonzero element of $V\otimes V$, $R$ is the one-dimensional subspace of $V\otimes V$ generated by $f$, $(x_1, \ldots, x_n)$ is a basis of $V$, $M$ is the matrix of coefficients of $f$ relative to this basis. Our aim is to examine whether a PBW deformation $U$ of $A=A(V,R)$ is Calabi-Yau if $A$ is Calabi-Yau (for Calabi-Yau algebras, the reader is referred to~\cite{g:cy} and to references therein). When the global dimension of $A$ is not 2, the PBW deformations of $A$ are described by Proposition \ref{infinite} with $N=2$, and $A$ is never Calabi-Yau because the global dimension of $A$ is infinite (\!\!~\cite{bt:cypbw}, Remark 2.8). 

In all the sequel, we assume that the global dimension of $A$ is 2, hence $n$ is $\geq 2$. Since $A$ is Calabi-Yau if and only if $A$ is the preprojective algebra of a non-Dynkin quiver (\!\!~\cite{boc:nccone}, Theorem 2.5), $A$ is Calabi-Yau if and only if $M$ is invertible and antisymmetric (see also Proposition \ref{Mant} below). For the moment, we do not impose such assumptions on $M$ and we want to describe the bimodule Koszul resolution and its dual for any $M$, the global dimension of $A$ still being equal to 2.

We need some notations. For every $a$, $b$, $c$ in $A$, set
$$c\stackrel{\ell}{.} (a\otimes b)=(a\otimes b)\stackrel{\ell}{.} c=a\otimes cb, \ \ 
c\stackrel{r}{.} (a\otimes b)=(a\otimes b)\stackrel{r}{.} c=ac\otimes b.$$
Then $A\otimes A$ is a left or right $A$-module for the action $\stackrel{\ell}{.}$ or $\stackrel{r}{.}$ respectively. In both cases, $A$ acts $A\otimes A$-linearly for the outer structure. In both cases, we need a notation where $c$ is on the left of $a\otimes b$ and another one where $c$ is on the right. 

Denote by $(A\otimes A)^{p\times q}$ the $A\otimes A$-module (for the outer structure) of the $p\times q$ matrices with entries in $A\otimes A$. Let $C=(c_{ij})$ be an $m\times p$ matrix with entries in $A$. The $A\otimes A$-linear maps for the outer structure
$$(A\otimes A)^{p\times q} \stackrel{C\stackrel{\ell}{.}}{\longrightarrow} (A\otimes A)^{m\times q}, \ (A\otimes A)^{p\times q} \stackrel{C\stackrel{r}{.}}{\longrightarrow} (A\otimes A)^{m\times q} $$
are defined by 
$$C\cdot \left( (u_{jk})_{1\leq j \leq p, 1\leq k \leq q} \right)= \left( \sum_{1\leq j \leq p} c_{ij} \cdot u_{jk}\right)_{1\leq i \leq m, 1\leq k \leq q},$$
where $\cdot$ denotes $\stackrel{\ell}{.}$ or $\stackrel{r}{.}$ respectively, and $(u_{jk})$ is a $p\times q$ matrix with entries in $A\otimes A$. The $A\otimes A$-linear maps
$$(A\otimes A)^{p\times q} \stackrel{\cdot C}{\longrightarrow} (A\otimes A)^{p\times r},$$
are defined similarly for any $q\times r$ matrix $C=(c_{jk})$ with entries in $A$ by
$$\left( (u_{ij})_{1\leq i \leq p, 1\leq j \leq q} \right) \cdot C= \left( \sum_{1\leq j \leq q} u_{ij} \cdot c_{jk}\right)_{1\leq i \leq p, 1\leq k \leq r}.$$

On the one hand, tensoring the left Koszul resolution (\ref{res2}) by $A$ on the right, we get
\begin{equation} \label{bir1}
0 \longrightarrow A\otimes A \stackrel{\stackrel{r}{.} x^tM}{\longrightarrow} (A\otimes A)^{1\times n} \stackrel{\stackrel{r}{.} x}{\longrightarrow} A\otimes A \longrightarrow 0,
\end{equation}
where $x$ is the $n\times 1$ matrix with entries $x_1, \ldots , x_n$ and $x^t$ is its transpose. Similarly, tensoring the right Koszul resolution by $A$ on the left, we get
\begin{equation} \label{bir2}
0 \longrightarrow A\otimes A \stackrel{Mx \stackrel{\ell}{.}}{\longrightarrow} (A\otimes A)^{n\times 1} \stackrel{x^t\stackrel{\ell}{.}}{\longrightarrow} A\otimes A \longrightarrow 0.
\end{equation}
In order to combine (\ref{bir1}) and (\ref{bir2}), we replace (\ref{bir2}) by its transpose
\begin{equation} \label{bir3}
0 \longrightarrow A\otimes A \stackrel{\stackrel{\ell}{.} x^tM^t}{\longrightarrow} (A\otimes A)^{1\times n} \stackrel{\stackrel{\ell}{.} x}{\longrightarrow} A\otimes A \longrightarrow 0.
\end{equation}
Following Proposition 3.1 in~\cite{vdb:hom} or Theorem 4.4 in~\cite{rbnm:kogo}, the bimodule Koszul resolution of $A$ is 
\begin{equation} \label{bir4}
\xymatrix@=1cm{0 \ar[r] & A\otimes A \ar[rr]^{\!\!\!\! \stackrel{r}{.} x^tM+\stackrel{\ell}{.} x^tM^t} &&(A\otimes A)^{1\times n}\ar[r]^{\ \ \ \ \stackrel{r}{.} x-\stackrel{\ell}{.} x} & A\otimes A \ar[r] & 0},  
\end{equation}
which is a minimal projective resolution of $A$ in the category $A$-grMod-$A$ of graded $A$-$A$-bimodules, for the multiplication $\mu : A\otimes A \rightarrow A$ of $A$.

On the other hand, tensoring the right dual Koszul complex (\ref{res3}) by $A$ on the left, we get
\begin{equation} \label{bir5}
0 \longrightarrow A\otimes A \stackrel{x \stackrel{\ell}{.}}{\longrightarrow} (A\otimes A)^{n\times 1} \stackrel{x^t M\stackrel{\ell}{.}}{\longrightarrow} A\otimes A \longrightarrow 0,
\end{equation}
which becomes after transposing it,
\begin{equation} \label{bir6}
0 \longrightarrow A\otimes A \stackrel{\stackrel{\ell}{.} x^t}{\longrightarrow} (A\otimes A)^{1\times n} \stackrel{\stackrel{\ell}{.} M^tx}{\longrightarrow} A\otimes A \longrightarrow 0.
\end{equation}
Combining the latter complex with the left dual Koszul complex (tensored by $A$ on the right), we obtain the \emph{dual bimodule Koszul complex}
\begin{equation} \label{bir8}
\xymatrix@=1cm{0 \ar[r] & A\otimes A \ar[r]^{\!\!\!\!\stackrel{r}{.} x^t-\stackrel{\ell}{.} x^t} &(A\otimes A)^{1\times n}\ar[rr]^{\ \ \ \ \stackrel{r}{.} Mx+\stackrel{\ell}{.} M^tx} && A\otimes A \ar[r] & 0}.  
\end{equation}

This complex in the category $A$-grMod-$A$ computes the Hochschild cohomology spaces HH$^i(A, A\otimes A)$ for $i=0,1,2$ (the other spaces are zero). In general, the map $\mu \circ (\stackrel{r}{.} Mx+\stackrel{\ell}{.} M^tx)$ is nonzero, unless $M$ is antisymmetric. Moreover, if $M$ is antisymmetric, we have the following commutative diagram
\begin{equation} \label{diag}
\xymatrix@=1cm{A\otimes A \ar[d]^{id}\ar[rr]^{\stackrel{r}{.} x^t-\stackrel{\ell}{.} x^t}&&(A\otimes A)^{1\times n}\ar[d]^{\cdot M}\ar[rr]^{\ \ \ \ \stackrel{r}{.} Mx-\stackrel{\ell}{.} Mx}&&A\otimes A\ar[d]^{id}\ar[r]^{\mu}& A \ar[d]^{id} \\
A\otimes A\ar[rr]^{\!\!\!\! \stackrel{r}{.} x^tM-\stackrel{\ell}{.} x^tM}&&(A\otimes A)^{1\times n}\ar[rr]^{\stackrel{r}{.} x-\stackrel{\ell}{.} x}&&A\otimes A \ar[r]^{\mu}& A }
\end{equation}
in which we have omitted the zero maps at the two ends for both rows. Since the entries of $M$ are in $k$, the map $\cdot M$ coincides with $\stackrel{r}{.}M$ and with $\stackrel{\ell}{.}M$.
\Bpo \label{Mant}
Assume $M$ is antisymmetric. Then \emph{HH}$^0(A, A\otimes A)=0$, and $A$ is $2$-Calabi-Yau if and only if $M$ is invertible.
\Epo
\Bdm
The first assertion is clear since the second row of (\ref{diag}) (including the zero maps) is exact. By definition, $A$ is 2-Calabi-Yau if and only if the first row (including the zero maps) is exact. Thus $A$ is 2-Calabi-Yau if $M$ is invertible. Conversely assume that $A$ is 2-Calabi-Yau. Using the same arguments as in the proof of Lemma \ref{lem}, we see that the first row is a minimal projective resolution in $A$-grMod-$A$ and thus $M$ is invertible.
\qed
\\

Let us keep our general assumptions: $M$ is arbitrary and the global dimension of $A$ is 2. Let $U$ be a PBW 
deformation of $A$. For the standard graded/filtered results implicitly used below, we refer to~\cite{nvo:grad}, Chapter D. For example, the global dimension of $U$ is $\leq 2$. According to Proposition \ref{gldim2}, we denote by $f-v-\lambda$ the relation of $U$ where $v\in V$ and $\lambda \in k$. Set $v=\sum_{1\leq i \leq n} \lambda_i x_i$ and let $\bar{v}$ be the $n \times 1$ matrix with entries $\lambda_1, \ldots , \lambda_n$. It is easy to check that the sequence of maps 
\begin{equation} \label{ubir}
\xymatrix@=1cm{0 \ar[r] & U\otimes U \ar[rr]^{\!\!\! \stackrel{r}{.} x^tM+\stackrel{\ell}{.} x^tM^t- \cdot \bar{v}^t} &&(U\otimes U)^{1\times n}\ar[r]^{\ \ \ \ \stackrel{r}{.} x-\stackrel{\ell}{.} x} & U\otimes U \ar[r]^{\mu}& U \ar[r] & 0}  
\end{equation}
where $\mu$ still denotes the multiplication, is a complex in the category $U$-filtMod-$U$ of filtered $U$-$U$-bimodules. Since the graded complex associated to (\ref{ubir}) is exact, (\ref{ubir}) is a projective resolution of $U$ in $U$-filtMod-$U$, called the \emph{bimodule Koszul resolution} of $U$.

For $\lambda =0$, the bimodule Koszul resolution of $U$ coincides with the standard complex defined by Braverman and Gaitsgory (\!\!~\cite{bg:pbw}, 5.4). For $n=2$ and $M$ antisymmetric, $U$ is a Sridharan enveloping algebra and the bimodule Koszul resolution of Sridharan enveloping algebras was described in~\cite{k:algenv}, Proposition 3. For $n>2$ or $M$ not antisymmetric, $A$ is not a polynomial algebra, hence $U$ is not a Sridharan enveloping algebra, even up to a linear change of generators.

The \emph{dual bimodule Koszul complex} of $U$ is the following
\begin{equation} \label{udual}
\xymatrix@=1cm{0 \ar[r] & U\otimes U \ar[r]^{\!\!\!\!\stackrel{r}{.} x^t-\stackrel{\ell}{.} x^t} &(U\otimes U)^{1\times n}\ar[rr]^{\ \ \ \stackrel{r}{.} Mx+\stackrel{\ell}{.} M^tx-\cdot \bar{v}} && U\otimes U \ar[r] & 0},  
\end{equation}
and computes the Hochschild cohomology spaces HH$^i(U, U\otimes U)$ for $i=0,1,2$ (the other spaces are zero). If $M$ is antisymmetric, Proposition \ref{Mant} shows that HH$^0(U, U\otimes U)=0$. In general the map $\mu \circ (\stackrel{r}{.} Mx+\stackrel{\ell}{.} M^tx-\cdot \bar{v})$ is nonzero, unless $M$ is antisymmetric and $v=0$. 

We are ready now to investigate further $U$ when $A$ is Calabi-Yau, i.e., when $M$ is antisymmetric and invertible. Then the graded complex associated to (\ref{udual}) is exact except perhaps at the last position on the right, so that we have 
$$\mathrm{HH}^0(U, U\otimes U)=\mathrm{HH}^1(U, U\otimes U)=0.$$
If $v=0$, we can replace $A$ by $U$ in the commutative diagram (\ref{diag}), thus $U$ is 2-Calabi-Yau in this case. Conversely, assume that $U$ is 2-Calabi-Yau. Then there exists a surjective $U\otimes U$-linear map $\mu _u : U\otimes U \rightarrow U$ with $\mu _u (1\otimes 1) =u \in U$, such that 
$$\xymatrix@=1cm{(U\otimes U)^{1\times n}\ar[rr]^{\ \ \ \stackrel{r}{.} Mx-\stackrel{\ell}{.} Mx-\cdot \bar{v}} && U\otimes U \ar[r]^{\ \ \ \mu _u}& U \ar[r] & 0}$$
is exact. Therefore $\sum_{1\leq i \leq n} \lambda_i a_i u b_i=0$ for any $a_i$ and $b_i$ in $U$. For each $i$, choosing $a_j=b_j=\delta_{ij}.1$ implies that $\lambda_i=0$ since $u\neq 0$. Thus $v=0$. We formulate what we have obtained as follows.
\Bte \label{cypbw}
Let $k$ be a field, let $n$ be an even natural number $\geq 2$, and let $U$ be the associative $k$-algebra defined by generators $x_1, \ldots , x_n$ subject to the single relation 
$$\sum_{1\leq i \leq n/2} [x_i, x_{i+n/2}] = v + \lambda,$$
where the bracket stands for the commutator, $v$ is a linear combination of the $x_i$'s, and $\lambda \in k$.
Then the filtered algebra $U$ is Koszul and we have $\mathrm{HH}^i(U, U\otimes U)=0$ whenever $i\neq 2$. Furthermore, $U$ is $2$-Calabi-Yau if and only if $v=0$.
\Ete

Let us note in particular that the PBW deformation $U$ obtained by adding only a constant to the relation of $A$ is Calabi-Yau. This phenomenon occurs in the following general context.
\Bte \label{cyconst}
Let $A=A(V,R)$ be an $N$-Koszul graded algebra, $N\geq 2$, with $V$ finite-dimensional. Let $\varphi : R \rightarrow k$ be $k$-linear (i.e., $\varphi= \varphi_0$). Assume that $U=T(V)/I(P)$ is a PBW deformation of $A$, where $P=(\mathrm{Id}-\varphi)(R)$. If $A$ is d-Calabi-Yau for a certain $d\geq 2$, then $U$ is $d$-Calabi-Yau. For example, any Sridharan enveloping algebra of an $n$-dimensional abelian Lie algebra is $n$-Calabi-Yau; in particular the Weyl algebra $A_n$ is $2n$-Calabi-Yau.
\Ete
\Bdm
It is known that $A$ is AS-Gorenstein of dimension $d$, $A$ satisfies the Van den Bergh duality and the automorphism of $A$ defining this duality is the identity automorphism (\!\!~\cite{bt:cypbw}, Section 4). Therefore we can use the intrinsic explicit isomorphism $\Phi'_{r-l}$ from the dual bimodule Koszul resolution $(A\otimes A^! \otimes A,\, \delta^\ast )$ to the bimodule Koszul resolution $(A\otimes A^{!\ast} \otimes A,\, \delta )$, as defined in \!\!~\cite{rbnm:kogo} (in the paragraph following Proposition 6.4). Denote $\Phi'_{r-l}$ by $\Phi_A$ for simplicity. Let us notice that $\Phi_A$ is uniquely determined by its restriction to $A^!$ and this restriction is mapped to $A^{!\ast}$. So the $U\otimes U$-linear map $\Phi_U : U\otimes A^! \otimes U \rightarrow U\otimes A^{!\ast} \otimes U$ coinciding with $\Phi_A$ on $A^!$ is well-defined. Using a graded/filtered argument, we see that $\Phi_U$ is an isomorphism of bimodules. 

Next, we recall that a bimodule Koszul resolution of $U$ was defined in~\cite{bg:hsra}, Section 5. This resolution is of the form $(U\otimes A^{!\ast} \otimes U,\, \delta )$ where the restriction of $\delta$ to $A^{!\ast}$ coincides with the restriction of the differential $\delta$ relative to $A$; the latter fact uses the equalities $U_i=F^i$ and $A_i=V^{\otimes i}$ whenever $i<N$. Consequently, the dual bimodule Koszul complex of $U$ is of the form $(U\otimes A^! \otimes U,\, \delta ^\ast)$ where the restriction of $\delta ^{\ast}$ to $A^!$ coincides with the restriction of the differential $\delta ^{\ast}$ relative to $A$. Thus, since $\Phi_A$ is a morphism of complexes, $\Phi_U$ is a morphism of complexes as well. The isomorphism of complexes $\Phi_U$ allows us to conclude that $U$ is $d$-Calabi-Yau.

We can take for $A$ a polynomial algebra over $n$ variables, so that $A$ is 2-Koszul and $n$-Calabi-Yau. Remark that $U$ is a PBW deformation for any $k$-linear $\varphi : R \rightarrow k$. In general, the Sridharan enveloping algebras are exactly the PBW deformations of polynomial algebras. Assuming that the ground field $k$ has zero characteristic, the linear part of such a deformation is an arbitrary Lie bracket, and $\varphi$ is an arbitrary 2-cocycle for this Lie bracket~\cite{k:algenv}. The particular case of the Weyl algebra $A_n$ corresponds to $2n$ variables and a symplectic form for $\varphi$.
\qed
\\

An interesting question is to find necessary and sufficient conditions on the Lie algebra and on the 2-cocycle $\varphi$ for an arbitrary Sridharan enveloping algebra to be Calabi-Yau. It would be also interesting to extend Theorem \ref{cypbw} to the preprojective algebra of a non-Dynkin quiver having more than one vertex. Other natural questions are the following: how can the results of Section 3 be extended to symmetric or supersymmetric $f$? Is Gerasimov's theorem true when the ground field $k$ is replaced by a semisimple ring? It is already known that Backelin's theorem (Proposition \ref{back}) holds for $k$ semisimple (\!\!~\cite{bgso:kdp}, Theorem 2.6.1, Remark).

\bigskip
\bigskip
\noindent
\emph{Roland Berger\\
Laboratoire de Math\'ematiques de l'Universit\'e de Saint-Etienne\\
Facult\'{e} des Sciences et Techniques\\ 23, Rue Docteur Paul Michelon, 42023 Saint-Etienne Cedex 2\\
France\\}
Roland.Berger@univ-st-etienne.fr


\begin{thebibliography}{99}
\bibitem{as:regular} M. Artin, W. F. Schelter, \emph{Graded algebras of global dimension 3}, Adv.
Math. \textbf{66} (1987), 171-216.
\bibitem{atv:modules} M. Artin, J. Tate, M. Van den Bergh, \emph{Modules over regular algebras of dimension 3}, Inv. Math. \textbf{106} (1991), 335-388.
\bibitem{bac:ratio} J. Backelin, \emph{La s\'erie de Poincar\'e-Betti d'une alg\`ebre gradu\'ee de type fini \`a une relation 
est rationnelle}, C. R. Acad. Sci. Paris, Ser A \textbf{287} (1978), 843-846.
\bibitem{bafr:kalg} J. Backelin, R. Fr\"{o}berg, \emph{Koszul algebras, Veronese subrings and rings with
linear resolutions}, Rev. Roum. Math. Pures Appli. \textbf{30} (1985), 85-97.
\bibitem{bgso:kdp} A. A. Beilinson, V. Ginzburg, W. Soergel, \emph{Koszul duality patterns in
representation theory}, J. Am. Math. Soc. \textbf{9} (1996), 473-527.
\bibitem{rb:conkos} R. Berger, \emph{Confluence and Koszulity}, J. Algebra
\textbf{201} (1998), 243-283.
\bibitem{rb:nonquad} R. Berger, \emph{Koszulity for nonquadratic algebras}, J. Algebra
\textbf{239} (2001), 705-734.
\bibitem{rbdvmw:homog} R. Berger, M. Dubois-Violette, M. Wambst, \emph{Homogeneous algebras},  J. Algebra \textbf{261} (2003), 172-185.
\bibitem{bg:hsra} R. Berger, V. Ginzburg, \emph{Higher symplectic reflection algebras and non-homogeneous N-Koszul property}, J. Algebra \textbf{304} (2006), 577-601.
\bibitem{rbnm:kogo} R. Berger, N. Marconnet, \emph{Koszul and Gorenstein properties for homogeneous algebras}, 
Algebr. Represent. Theory \textbf{9} (2006), 67-97.
\bibitem{bt:cypbw} R. Berger, R. Taillefer, \emph{Poincar\'e-Birkhoff-Witt deformations of Calabi-Yau algebras}, J. Noncommut. Geom. \textbf{1} (2007), 241-270.
\bibitem{bock:gcy} R. Bocklandt, \emph{Graded Calabi Yau algebras of dimension 3} 
(with an appendix by M. Van den Bergh), J. Pure Appl. Algebra  \textbf{212} (2008), 14-32.
\bibitem{boc:nccone}  R. Bocklandt, \emph{Noncommutative Tangent Cones and Calabi Yau Algebras}, arXiv:0711.0179.
\bibitem{bou:alghom} N. Bourbaki, \emph{Alg\`ebre, Chapitre 10}, Masson, Paris, 1980.
\bibitem{bg:pbw} A. Braverman, D. Gaitsgory, \emph{Poincar\'e-Birkhoff-Witt theorem for quadratic algebras of Koszul type}, J. Algebra \textbf{181} (1996), 315-328.
\bibitem{bbk:periodic} S. Brenner, M. C. R. Butler, A. D. King,\emph{Periodic algebras which are almost Koszul},  Algebr. Represent. Theory \textbf{5} (2002), 331-367.
\bibitem{cw:bili} B. Corbas, G. D. Williams, \emph{Bilinear forms over an algebraically closed field}, J. Pure Appl. Algebra \textbf{165} (2001), 255-266.
\bibitem{mdv:multi} M. Dubois-Violette, \emph{Multilinear forms and graded algebras}, J. Algebra \textbf{317} (2007), 198-225.
\bibitem{dvp:plac} M. Dubois-Violette, T. Popov, \emph{Homogeneous algebras,
  statistics and combinatorics}, Lett. Math. Phys. \textbf{61} (2002), 159-170.
\bibitem{ep:algMMT} P. Etingof, I. Pak, \emph{An algebraic extension of the MacMahon master theorem}, Proc. Amer. Math. Soc. \textbf{136} (2008), 2279-2288.
\bibitem{fv:pbwdef} G. Fl\o ystad, J. E. Vatne, \emph{PBW-deformations of N-Koszul algebras},
J. Algebra \textbf{302} (2006), 116-155.
\bibitem{fro:mono} R. Fr\"{o}berg, \emph{Determination of a class of Poincar\'{e} series}, Math.
Scand. \textbf{37} (1975), 29-39.
\bibitem{gera:distrib} V. N. Gerasimov, \emph{Free associative algebras and inverting homomorphisms of rings}, 
Amer. Math. Soc. Transl. \textbf{156} (1993).
\bibitem{g:cy} V. Ginzburg, \emph{Calabi-Yau algebras}, arXiv:math.AG/0612139.
\bibitem{gm:delta} E. L. Green, E. N. Marcos, \emph{$\delta$-Koszul algebras}, Comm. Algebra \textbf{33} (2005), 1753-1764. 
\bibitem{gmmz:dkos} E. L. Green, E. N. Marcos, R. Mart\'{\i}nez-Villa, Pu Zhang, \emph{D-Koszul algebras}, J. Pure Appl. Algebra \textbf{193} (2004), 141-162. 
\bibitem{gs:fingen} E. L. Green, N. Snashall, \emph{Finite generation of Ext for a generalization of $D$-Koszul algebras}, J. Algebra \textbf{295} (2006), 458-472.  
\bibitem{phhml:NMMT} P. H. Hai, B. Kriegk, M. Lorenz, \emph{N-homogeneous superalgebras}, J. 
Noncommut. Geom. \textbf{2} (2008), 1-51.
\bibitem{phhml:2MMT} P. H. Hai, M. Lorenz, \emph{Koszul algebras and the quantum MacMahon master theorem}, 
Bull. London Math. Soc. \textbf{39} (2007), 667-676.
\bibitem{k:algenv} C. Kassel, \emph{L'homologie cyclique des alg\`ebres enveloppantes}, Invent. math. \textbf{91} (1988), 221-251.
\bibitem{kriegk:cras} B. Kriegk, \emph{Un crit\`ere num\'erique pour la propri\'et\'e de Koszul g\'en\'eralis\'ee}, 
C. R. Acad. Sci. Paris, Ser I \textbf{344} (2007), 545-548.
\bibitem{nvo:grad} C. Nastacescu, F. Van Oystaeyen, \emph{Graded ring theory},
North-Holland, 1982.
\bibitem{pp:quad} A. Polishchuk, L. Positselski, \emph{Quadratic algebras}. Univ. Lecture Ser., vol. 37. 
Amer. Math. Soc., Providence, RI, 2005.
\bibitem{vdb:hom} M. Van den Bergh, \emph{Non-commutative homology of some three-dimensional quantum
spaces}, K-Theory \textbf{8} (1994), 213-220.
\end{thebibliography}
\end{document}